\documentclass[11pt]{article}
\usepackage{amssymb}
\usepackage{latexsym}
\input amssym.def
\input amssym
\newtheorem{theorem}{Theorem}[section]
\newtheorem{lemma}[theorem]{Lemma}

\newtheorem{proposition}[theorem]{Proposition}
\newtheorem{corollary}[theorem]{Corollary}

\newtheorem{remark}{Remark}

\newcommand{\ds}{\displaystyle}

\newcommand{\halmos}{\rule{1ex}{1.4ex}}

\newcommand{\bea}{\begin{eqnarray}}
\newcommand{\epf}{\hspace*{\fill}\mbox{$\halmos$}}
\newcommand{\eea}{\end{eqnarray}}
\newcommand{\nn}{\nonumber \\}
\newcommand{\be}{\begin{equation}}

\newcommand{\ee}{\end{equation}}

\usepackage{amssymb}

\title{{Virasoro Algebra, Dedekind $\eta$-function and Specialized Macdonald's Identities} }
\author{Antun  Milas}

\begin{document}
\date{}
\bibliographystyle{alpha}
\maketitle
\begin{abstract}
\noindent 
We motivate and prove a series of identities which form a  
generalization of the Euler's pentagonal number theorem, and are
closely related to specialized Macdonald's identities for powers of the Dedekind $\eta$--function. More precisely, we show that what we call ``denominator formula'' for the
Virasoro algebra has ``higher analogue'' for 
all $c_{s,t}$-minimal models. We obtain one identity per series which is in 
agreement with features of conformal field theory such as {\em fusion} and 
{\em modular invariance} that require 
all the irreducible modules of the series. 
In particular, in the case of $c_{2,2k+1}$--minimal 
models we give a new proof of a family of specialized Macdonald's identities associated with twisted affine Lie algebras of type $A^{(2)}_{2k}, k \geq 2$ (i.e., $BC_k$-affine root system) 
which involve $(2k^2-k)$-th powers of the Dedekind $\eta$-function.
Our paper is in many ways a continuation of math.QA/0309201.
\end{abstract}
\renewcommand{\theequation}{\thesection.\arabic{equation}}
\renewcommand{\thetheorem}{\thesection.\arabic{theorem}}
\setcounter{equation}{0} \setcounter{theorem}{0}
\setcounter{section}{0}

\section{Introduction} In 1972 I.G. Macdonald discovered a beautiful series of 
identities associated to affine root systems \cite{Ma}.
From the affine Kac--Moody algebra point of view, Macdonald's
identities are well understood and are instances of 
what is nowadays called {\em
denominator identity} \cite{Ka}, \cite{Ka1}, often written in the
following form: \be \label{mac} \sum_{w \in W} \epsilon(w)
e^{w(\rho)-\rho}=\prod_{ \alpha \in \Delta_+}
(1-e^{-\alpha})^{{\rm mult}(\alpha)}, \ee where $\Delta$ is an
affine root system, $\epsilon(w) \in \{\pm 1\}$, $\Delta_+$ is the set of  positive roots, $W$
is the affine Weyl group, ${\rm
mult}(\alpha)$ are the root multiplicities and $\rho$ denotes the half
sum of the positive roots \cite{Ka1}. 
By using various
specializations Macdonald obtained several $q$-identities, the
most interesting of which involve
$$\eta(q)^{{\rm dim}(\goth{g})}$$
$$\eta(q)=q^{1/24} \prod_{i \geq 1} (1-q^i)$$
(i.e., Dyson-Macdonald's identities \cite{Dy}), 
appearing on the right hand side of (\ref{mac}), where ${\rm
dim}(\goth{g})$ denotes the dimension of a finite dimensional
simple Lie algebra $\goth{g}$, such as for example
\be \label{list}
{\rm dim}(\goth{g})=3, 8, 10, 14, 15, 21, 24,...
\ee
For instance, in the case of $\goth{g}=\goth{sl}_2$, a specialized
Macdonald's identity gives a famous Jacobi's identity
$$\eta(q)^3=q^{1/8}\sum_{m=0}^\infty (-1)^m
(2m+1)q^{\frac{m(m+1)}{2}}.$$

By now there are several different proofs, interpretations, reinterpretations and
extensions of Macdonalds's work. 
Let us mention a few most important contributions.
A new hat for these identities was obtained 
by Kac \cite{Ka} who placed Macdonald's identities in the
context of affine Lie algebras (cf. \cite{Mo}) as a very special case of the
Weyl-Kac character formula for the trivial module (cf. \cite{Ka1}). Garland and Lepowsky generalized
an earlier work by Kostant to general Kac-Moody Lie algebras \cite{GL} (cf. \cite{G}), which, in particular, led to Macdonald's identities. More recently, some 
methods and results from \cite{GL} were 
generalized by Borcherds so that they apply to generalized
Kac--Moody algebras (for an important application of this theory
see \cite{Bo}). There are many different approaches and proofs of Macdonald's identities 
and generalizations, such as \cite{As}, \cite{Fe}, \cite{Fe1}, \cite{Ko}, \cite{Lo}, \cite{Ka2}, \cite{Le0}, \cite{Le}, \cite{Fr}, etc. 
In particular, the main result in \cite{Le0} indicates
that distinguished powers (\ref{list}) are by no means distinguished from
the affine Kac--Moody algebra point of view (see also \cite{Ka1}). 
(We thank Jim Lepowsky for pointing us to \cite{Le0}.) 

Unlike the previous approaches, in this paper we obtain some 
specialized Macdonald's identities by using the 
representation theory of the Virasoro algebra.
Compared to affine Lie algebras \cite{PS}, the Virasoro algebra
has a simple geometric interpretation; being the 
non-trivial central extension of the Lie algebra of polynomial
vector fields on the circle, $Vect(S^1)$. 
In spite of this simplicity the highest (or the lowest) weight representation theory turns out to be quite
complicated and the most interesting representations of the Virasoro algebra have no
(known) geometric realization. Nevertheless, a complete
classification of irreducible highest weight modules (including
BGG-type resolutions) was given by Feigin and Fuchs \cite{FFu1},
\cite{FFu} and Rocha-Caridi and Wallach who obtained some partial
results prior to Feigin and Fuchs' results \cite{RCW}, \cite{RCW1}
($c=1$ case is due to Kac \cite{KR}). Actually, Rocha-Caridi was
the first \cite{RC} who wrote explicitly formulas for graded
dimensions (or characters) of irreducible highest weight Virasoro algebra modules
(after \cite{FFu1}). Among all the highest weight irreducible modules the most interesting are those 
parameterized by the {\em central charge} (cf. \cite{KR})
\be \label{cc}
c_{s,t}=1-\frac{6(s-t)^2}{st}, \  s,t \in \mathbb{N}, \ s,t \geq 2, \ (s,t)=1,
\ee
i.e., {\em minimal models}. 
One of the main reasons why the Virasoro algebra and minimal models have been
studied intensively over the last two decades is because of their
relevance in {\em conformal field theory} \cite{BPZ}, \cite{S} (but also in
2-dimensional statistical physics and integrable models). In addition, the
Virasoro algebra is related to affine Lie algebras via
Sugawara construction (see for instance \cite{Ka1},\cite{KR}),
which allows us to obtain all the unitary minimal
models (cf. \cite{KR}) via {\em coset} constructions \cite{GKO},
\cite{KW1} from the integrable highest weight modules for
affine Lie algebras. For non-unitary minimal models (considered for instance in this paper) there is a version of
GKO-construction which uses cosets associated to representations
of affine Lie algebras at admissible levels (certain rational 
levels above the critical level \cite{KW}). These
admissible level modules are rather mysterious so it is not clear
how the work \cite{KW} helps in understanding non-unitary minimal
models. Anyhow, it is known that the
minimal models with the central charge $c_{s,t}$ are sources of 
{\em rational vertex operator algebras} \cite{Wa}, \cite{Zh} (even better, these vertex operator algebras
are {\em regular} \cite{DLM}), genus-zero weakly holomorphic 
conformal field theories \cite{Hu0}, modular invariant theories \cite{Zh}, \cite{Hu}, even modular
functors \cite{BFM}, so clearly these objects are of crucial importance
(cf. \cite{Fel}, \cite{BMCS} for different approaches to Virasoro minimal models)

At the abstract level, as in the affine Lie algebra case, the Virasoro
algebra has what we call {\em denominator formula}; a consequence of 
the Euler-Poincar\'e principle applied to a resolution of the {\em trivial} Virasoro algebra module (with
$c=0$) in terms of Verma modules \cite{RCW}. This resolution gives a $q$--series identity 
equivalent 
to a classical Euler's formula (cf. \cite{A}) \be
\label{euler} \eta(q)=q^{1/24} \sum_{n \in \mathbb{Z}}
(-1)^n q^{\frac{3n^2-n}{2}}. \ee 
Notice that $c=0$ occurs on the list
(\ref{cc}) for $s=2,t=3$, and that this is the only minimal model with
$c=0$ (cf. Section 3).

In general, for every (irreducible) minimal model we can compute 
its graded trace or the {\em character}, but
a single irreducible module does not carry a full information unless
it ``interacts'' with other minimal models with the same central charge 
(e.g., {\em fusion} \cite{BPZ},\cite{FHL}, etc.). Therefore, 
as in the central charge zero case, we would like to have some conformal field
theoretical formula that takes into account {\em all} the minimal
models with the same central charge, generalizing the {denominator} formula (\ref{euler}). 

In the present paper we address the following question:
\vskip 3mm 
\noindent 
{\em What are the conformal field theoretical analogue of (\ref{euler}) for 
other $c_{s,t}$--series?} 
\vskip 3mm 

In Theorem \ref{theoremmain2} we give the answer for all
$c_{s,t}$--minimal models. However, in this paper 
we will be primarily interested 
in $c_{2,2k+1}$-minimal models. These models are important because of their combinatorial 
interpretation (e.g., Andrews-Gordon identities \cite{FFr} and dilogarithm identities \cite{FS}). 
For these series our ``denominator formulas'' are equivalent 
to a series of specialized Macdonald's identities associated with twisted affine Lie algebras 
of type $A_{2k}^{(2)}$, $k \geq 2$ (i.e., $BC_k$ affine root system \cite{Ma}).
\noindent We will prove the following theorem (essentially a formula on p.138, \cite{Ma}):
\begin{theorem} ($c_{2,2k+1}$-denominator formula) \label{theoremmain}
For every  $k \in \mathbb{N}$, $k \geq 2$ \be \label{bck}
\eta(q)^{2k^2-k}=C_{k} (-1)^{\frac{k(k-1)}{2}}
\sum_{{\bf n} \in \mathbb{Z}^k} (-1)^{\sum_{i=1}^k n_i}
\chi_D({\bf n}) q^{L({\bf n})}, \ee where 
${\bf n}=(n_1,...,n_k) \in \mathbb{Z}^k$, 
\be \label{bfn} L({\bf
n})=\frac{2k^2-k}{24}+\sum_{i=1}^k \left( \frac{(2k+1)n_i^2}{2}
+\frac{(2i-1)n_i}{2}\right), \ee
 \be \label{chior}
 \chi_D({\bf n})=\prod_{1 \leq i
<j \leq k} \left\{(2i-1+n_i(4k+2))^2-(2j-1+n_j(4k+2))^2 \right\} \ee
and
$$C_k=\frac{1}{2^{k(k-1)} \prod_{1 \leq i<j \leq k}(i-j)(i+j-1)}.$$
\end{theorem}
Notice that $k=1$ case is not included in our theorem, but
our motivation clearly indicates that the Euler's identity should be added 
at the beginning of this list of identities.

We also obtain a generalization of (\ref{theoremmain}) (see
Theorem \ref{theoremmain2}). However, at this point we do not
fully understand whether some of identities from Theorem \ref{theoremmain2} can 
be related to other specialized Macdonald's identities.

Let us elaborate our proof of Theorem \ref{theoremmain}.
\begin{itemize}
\item We apply vertex operator algebra methods to study
the characters of minimal models \cite{Zh}, \cite{M}, \cite{M1}, \cite{DMN} and derive 
differential equations with fundamental system of solutions formed
by characters of irreducible modules with 
$c=c_{2,2k+1}$, $k \geq 2$. 
\item By using
the Abel's lemma \ref{abel}, as in \cite{M}, we obtain a list of
identities that involve 
$$\eta(q)^{2k(k-1)}, \ \ k \geq 2.$$
\item We ``factor'' {\em missing} powers
of $\eta(\tau)$ from the denominators of characters of irreducible modules which
gives a series of identities for
$$\eta(q)^{(2k-1)k}, \ k \geq 2.$$
Notice that 
$${\rm dim}(D_k)={\rm dim} (\goth{so}(2k))=(2k-1)k, \ k \geq 1.$$
In the last step we evaluate the Wronskian and expresse the result in terms
of Vandermonde determinants. Surprisingly, our identities are not
associated to $D_k^{(1)}$--series but rather to $A_{2k}^{(2)}$--series.
\end{itemize}
In what follows: 
$\mathbb{H}$ is the upper half-plane, $q=e^{2 \pi i \tau}, \tau \in \mathbb{H}$, $\mathbb{N}$
is the set of positive integer and $\mathbb{N}_0$ is the set of non-negative integers. 

\renewcommand{\theequation}{\thesection.\arabic{equation}}
\renewcommand{\thetheorem}{\thesection.\arabic{theorem}}
\setcounter{equation}{0} \setcounter{theorem}{0}

\section{Determinants}
Let
$$ V(x_1,...,x_k)=\left| \begin{array}{ccccc} 1 & 1 & . & . & 1 \\
                              x_1 & x_2 & . & . & x_k \\
                              x^2_1 & x^2_2& . & . &  x^2_k \\
                              . & . & . & . &  . \\
                              x^{k-1}_1 & x^{k-1}_2 & . & . & x^{k-1}_k
\end{array} \right|$$
denote the Vandermonde determinant associated to $x_1,...,x_k$.
This determinant can be computed by using the well-known formula
$$V(x_1,...,x_k)=\prod_{1 \leq  j < i \leq k}
(x_i-x_j),$$ which is equivalent
to the Weyl denominator formula for the finite dimensional Lie algebra $\goth{sl}_{k}$.

In complex analysis another determinant plays a
prominent rule. Wronskian (or Wronski) determinant associated to a set of
analytic functions $y_1(\tau),...,y_k(\tau)$ is given by
$$W(y_1,...,y_k)=\left| \begin{array}{ccccc} y_1(\tau) & y_2(\tau) & . & . & y_k(\tau) \\
                              y'_1(\tau) & y'_2(\tau) & . & . & y'_k(\tau) \\
                              y''_1(\tau) & y''_2(\tau) & . & . &  y''_k(\tau) \\
                              . & . & . & . &  . \\
                              y^{(k-1)}_1(\tau) & y^{(k-1)}_2(\tau) & . & . & y^{(k-1)}_k(\tau)
\end{array} \right|$$
and it is important in the following fundamental result
due to Abel:
\begin{lemma} \label{abel}
Suppose that $f_i(\tau)$ are holomorphic functions in $\mathbb{H}$, and for every $k \geq 1$, let
$$y_1(\tau),...,y_k(\tau)$$
form a fundamental system of solutions for
$$\left(\frac{d}{d \tau} \right)^{k} y(\tau)+f_1(\tau)\left(\frac{d}{d \tau} \right)^{k-1}y(\tau)+
\cdots+f_k(\tau) y(\tau)=0,$$
then
$$W(y_1,...,y_k)=Ce^{\ds{-\int^\tau_{\tau_0} f_1(\tau) d \tau}},$$
where $C=W(y_1(\tau_0),...,y_k(\tau_0))$ and  $\tau_0 \in \mathbb{H}$ (in fact, $C$ is some 
non-zero constant which depends on $y_1,...,y_k$ and $\tau_0$,
but not on $\tau$).
\end{lemma}
Of course, as in the case of differential equations, the Wronskian
can be used to determine whether the set of functions are linearly
independent.

In all our applications $y_i(\tau)$'s are analytic in the upper
half-plane and meromorphic at infinity (i.e., have $q$--expansions). 
Also, from now on
$$'=\left(\frac{1}{2\pi i} \frac{d}{d \tau}\right)=\left(q \frac{d}{dq} \right).$$
The next result shows that Vandermonde and Wronskian determinants
are closely related.
\begin{lemma}
Suppose that $y_i(\tau)$, $i=1,...,k$, are holomorphic in $\mathbb{H}$, with the $q$-expansions
$$y_i(q)=\sum_{n_i \geq \nu_i} a^{(i)}_{n_i} q^{n_i},$$
where $a^{(i)}_{n_i} \in \mathbb{C}$, $\nu_i \in \mathbb{Q}$, for
$i=1,...,k$. Then the Wronskian $W(y_1(\tau),...,y_n(\tau))$ is
holomorphic in $\mathbb{H}$ and its $q$-expansion at infinity is
given by \be \label{wronskivander} W(y_1,...,y_n)=\sum_{n_1 \geq
\nu_1,...,n_k \geq \nu_k} V(n_1,...,n_k) \left( \prod_{i=1}^k
a^{(i)}_{n_i} \right) q^{n_1+\cdots+n_k}. \ee
\end{lemma}
{\em Proof:} From the definition of Wronskian it is clear that
$W(y_1(\tau),..,y_k(\tau))$ is holomorphic and meromorphic at
infinity (it has a $q$-expansion). Thus, the only thing we have to show
is (\ref{wronskivander}). Let $A=\{a_{i,j}\}$ be a matrix of order
$k$ and $r \in \mathbb{N}$, $1 \leq r \leq k$. Suppose that
$a_{j,r}=b_{j,r}+c_{j,r}$ for every $j=1,...,k$, then
$${\rm det}(A)={\rm det}(B_r)+{\rm det}(C_r),$$
where $B_r$ and $C_r$ are matrices obtained from the matrix A by
replacing the $r$-th column by $[b_{1,r},...,b_{k,r}]$ and
$[c_{1,r},...,c_{k,r}]$, respectively. In $W(y_1,...,y_n)$ all our
entries are sums so we can repeat the previous formula for all the
columns and simultaneously factor $q^{n_i}$ from the $i$-th column
for every $i$. The remaining coefficient of $q^{n_1+\cdots+n_k}$ 
is the Vanderomonde
determinant for $n_1,...,n_k$. \epf
\renewcommand{\theequation}{\thesection.\arabic{equation}}
\renewcommand{\thetheorem}{\thesection.\arabic{theorem}}
\setcounter{equation}{0} \setcounter{theorem}{0}

\section{The Virasoro algebra}

In this part we recall a few basic results regarding the Virasoro algebra. 
For a good introduction to infinite-dimensional Lie algebra theory, the Virasoro
algebra and its minimal models see \cite{KR}. Let us recall that
$c_{s,t}$--minimal models are parameterized by the central charge
$$c_{s,t}=1-\frac{6(s-t)^2}{st},$$
where $s, t \in \mathbb{N}$, $ s,t \geq 2$, $(s,t)=1$, and by the
weights
$$h^{m,n}_{s,t}=\frac{(ns-mt)^2-(s-t)^2}{4 st},$$
where $1 \leq m <s$, $1 \leq n <t$. We also let
$$\bar{h}^{m,n}_{s,t}=h^{m,n}_{s,t}-\frac{c_{s,t}}{24}.$$
Once we fix $c_{s,t}$ there are exactly
$$\frac{(s-1)(t-1)}{2}$$
different values of $h_{s,t}^{m,n}$ for $1 \leq m <s$, $1 \leq n <t$. 
As in \cite{M}, we denote by $L(c,h)$ the irreducible highest
weight module with the central charge $c$ and the weight $h$, and
by
$$\bar{{\rm ch}}_{c,h}(q)={\rm tr}|_{L(c,h)} q^{L(0)-c/24}$$
the {\em graded dimension} or simply the {\em character} of
$L(c,h)$. The following result is well known (cf. \cite{RC})
\begin{theorem} \label{struct}
We have 
\be
\bar{{\rm ch}}_{c_{s,t}, h_{s,t}^{m,n}}(q)=\frac{q^{(h_{s,t}^{m,n}-c_{s,t}/24)} 
\ds{\sum_{r \in \mathbb{Z}}} \left( q^{st r^2+r(ns-mt)}-q^{str^2+r(ns+mt)+mn}\right)}{(q)_\infty}.
\ee
\end{theorem}
\epf

It is not hard to see that the expression
$$\sum_{r \in \mathbb{Z}} \left( q^{st r^2+r(ns-mt)}-q^{str^2+r(ns+mt)+mn} \right)$$
involves only non-negative powers of $q$.

To connect our results with \cite{M}, we shall first consider minimal models with
$$c_{2,2k+1}=1-\frac{6(2k-1)^2}{(4k+2)}, \ \ k \geq 2,$$
and the corresponding weights
$$h^{1,i}_{2,2k+1}=\frac{(2(k-i)+1))^2-(2k-1)^2}{8(2k+1)}, \ \ i=1,...,k.$$
The previous theorem implies a well-known result: For $i=1,...,k$, \be
\label{character} \bar{{\rm ch}}_{c_{2,2k+1},h^{1,i}_{2,2k+1}}(q)=
\frac{\ds{q^{(h^{1,i}_{2,2k+1}-\frac{c_{2,2k+1}}{24}+\frac{1}{24})}
\sum_{n \in \mathbb{Z}} (-1)^n
q^{\frac{(2(k-i)+1)n+(2k+1)n^2}{2}}}}{\eta(q)}. \ee In \cite{M} we
have used different, but equivalent formula (cf. \cite{FFr})
$$\bar{{\rm ch}}_{c_{2,2k+1},h^{1,i}_{2,2k+1}}(q)
=q^{(h^{1,i}_{2,2k+1}-\frac{c_{2,2k+1}}{24})} \prod_{n  \neq \pm i,0
\ {\rm mod}\ (2k+1)} \frac{1}{(1-q^{n})},$$ obtained from
(\ref{character}) by application of the Jacobi triple product
identity \cite{A}. We have used the infinite-product expressions (\ref{character}) 
in \cite{M} for purposes of proving one of Ramanujan's "Lost Notebook" formulas.
\renewcommand{\theequation}{\thesection.\arabic{equation}}
\renewcommand{\thetheorem}{\thesection.\arabic{theorem}}
\setcounter{equation}{0} \setcounter{theorem}{0}

\section{Modular invariance}
An interesting fact about minimal models is that, once we fix the
level $c_{s,t}$, the vector space spanned by graded traces
$$\{ \bar{ {\rm ch} }_{ c_{s,t},h^{m,n}_{s,t}}(q) \  : \ 1 \leq m <s, 1 \leq n <t  \}$$
is modular invariant \cite{CIZ}, i.e., a $SL(2,\mathbb{Z})$--module, where an
element $\gamma \in SL(2, \mathbb{Z})$ acts on the modulus $\tau$ in the standard way.
The best explanation of this phenomena was provided by Zhu in his work on modular invariance
of characters \cite{Zh} (see also \cite{Hu}), which uses in 
an essential way the theory of vertex operator algebras
\cite{FHL}, \cite{FLM}, \cite{Hu12}.
For purposes of this paper we do
not recall any of the theory of vertex operator algebras here, but we
mention that large portions of this theory have been used in \cite{M}
and henceforth in this paper.
\renewcommand{\theequation}{\thesection.\arabic{equation}}
\renewcommand{\thetheorem}{\thesection.\arabic{theorem}}
\setcounter{equation}{0} \setcounter{theorem}{0}

\section{``Strange formulas"}
We showed in \cite{M} that {\em product} of all characters of
$c_{2,2k+1}$--minimal models can be expressed in terms of
powers of quotients of two Dedekind eta functions with different moduli. 
One
of the reasons for multiplying these $q$--series stems from the
following observation (cf. \cite{M}): 
\be \label{etak}
\sum_{i=1}^k \bar{h}^{1,i}_{2,2k+1}=\sum_{i=1}^k
\left(h^{1,i}_{2,2k+1}-\frac{c_{2,2k+1}}{24} \right)
=\frac{2k(k-1)}{24}. \ee Here the number $k$ is also the 
number of inequivalent $c_{2,2k+1}$-minimal models (or the number of inequivalent
irreducible modules for the vertex operator algebra $L(c_{2,2k+1},0)$ \cite{Wa},
cf. \cite{M}). The rational number appearing on the right hand side
of (\ref{etak}) is important because it is related to asymptotic
behavior of products of characters as $q \rightarrow 0$.
Now, let us compute a version of (\ref{etak}) for arbitrary minimal models. 
As we already mentioned there are in
total
$$a_{s,t}=\frac{(s-1)(t-1)}{2}$$
inequivalent minimal models.
\noindent It is easy to see that  \bea \label{etakmore} &&
\frac{1}{2} \sum_{m=1}^{s-1} \sum_{n=1}^{t-1}
\bar{h}_{s,t}^{m,n}=\frac{1}{2} \sum_{m=1}^{s-1} \sum_{n=1}^{t-1}
\left(
\frac{(ns-mt)^2-(s-t)^2}{4st}-\frac{1}{24}\left(1-\frac{6(s-t)^2}{st}\right)\right)
\nn && =\frac{(s-1)(t-1)(st-s-t-1)}{48}, \eea where we rescale the
sum (viz. the factor $1/2$) because every number $h_{s,t}^{m,n}$ appears in the summation exactly twice.
Because of the identity
$$\frac{(s-1)(t-1)(st-s-t-1)}{48}=\frac{2 a_{s,t}(a_{s,t}-1)}{24},$$
it follows that (\ref{etakmore}) depends, as in the $c_{2,2k+1}$
case, only on the number of inequivalent minimal models.
Notice that the list $2k(k-1)$, $k \geq 2$ does not correspond to any of 
the lists of dimensions of classical finite-dimensional Lie algebras (i.e., 
${\rm dim}(A_k)$, ${\rm dim}(B_k)$, ${\rm dim}(C_k)$ and ${\rm dim}(D_k)$). Nevertheless,
if we add $\frac{1}{24}$ contribution from
each of the characters (\ref{character}), we get a much nicer expression 
\be \label{adding}
\frac{2k(k-1)+k}{24}=\frac{(2k-1)k}{24}, \ee 
which is equal to
$$\frac{{\rm dim}(\goth{so}(2k))}{24}$$
and appears on the Dyson-Macdonald's list.
The last expression indicates that there might be some kind of relationship between 
{\em higher dimensional} conformal field theory and the Virasoro 
minimal models.

In the next section
the observation (\ref{adding}) will be used in connection with
some specialized Macdonald's identities.
\renewcommand{\theequation}{\thesection.\arabic{equation}}
\renewcommand{\thetheorem}{\thesection.\arabic{theorem}}
\setcounter{equation}{0} \setcounter{theorem}{0}

\section{The Main Theorem}
In \cite{M}, for every integer $k$ we obtained a homogeneous differential 
equation of order $k$ with a fundamental system of solutions
formed by the characters of $c_{2,2k+1}$--minimal models. In our
approach the crucial role will play the following holomorphic,
quasimodular (normalized) Eisenstein series
$$\tilde{G}_{2}(q)=\frac{-1}{12}+2 \sum_{n \geq 1} \frac{n
q^n}{1-q^n}.$$
Let us recall a result from \cite{M} (see the proof of Theorem 8.6):
\begin{theorem} \label{diffmain}
There is a homogeneous differential equation of order $k$ with
holomorphic coefficients \be \label{diffm1} \left(q \frac{d}{dq}
\right)^k {y}(\tau)+k(k-1)\tilde{G}_2(\tau)
\left(q\frac{d}{dq}\right)^{k-1} {y}(\tau)+ \cdots =0, \ee with a
fundamental system of solutions formed by
$$y_i(\tau)=\bar{ch}_{c_{2,2k+1},h_{2,2k+1}^{1,i}}(q), \ i=1,...,k.$$
Moreover,
$$W(y_1,...,y_k)=\lambda_k \eta(\tau)^{2k(k-1)},$$
where $\lambda_k$ is a non-zero constant which depends only on $k$.
\end{theorem}
As we observed in the previous section, in order to achieve the
right powers of the Dedekind $\eta$--function as in Dyson-Macdonald's
identities for $D_k$--series we need additional powers of the Dedekind $\eta$--function.
These powers can be obtained by clearing out the denominator 
$${\eta(\tau)}$$
in the character formula
(\ref{character}). The resulting expressions (i.e., the numerators in (\ref{character})) 
satisfy another linear differential equation. 
The following lemma incorporates this procedure:
\begin{lemma} \label{fixing} 
After the substitution  $\tilde{y}(\tau)=y(\tau) \eta(\tau)$, the
homogeneous differential equation (\ref{diffm1}) becomes \be
\label{diffmod} \left(q \frac{d}{dq} \right)^k
\tilde{y}(\tau)+\left(k(k-1)+\frac{k}{2}\right)\tilde{G}_2(\tau)
\left(q\frac{d}{dq}\right)^{k-1}\tilde{y}(\tau)+ \cdots =0, \ee
with a fundamental system of solutions formed by
$$\tilde{y}_i(\tau), \ \ i=1,...,k.$$
\end{lemma}
{\em Proof:} Our starting point is the homogeneous differential
equation of degree $k$ in Theorem \ref{diffmain}.
By using the logarithmic derivative formula for the Dedekind eta
function \cite{A}, \cite{M}, we easily compute
$$\left(q \frac{d}{dq}\right) \tilde{y}(\tau)=
\left(q \frac{d}{dq}\right) \left( \eta(\tau) {y}(\tau)\right)$$
$$=-\frac{1}{2} \tilde{G}_2(\tau) \eta(\tau){y}(\tau)+\eta(\tau)\left(q \frac{d}{dq}\right)
y(\tau),$$ hence
$$\eta(\tau)\left(q \frac{d}{dq}\right) y(\tau)=\left(\left(q
\frac{d}{dq}\right)+\frac{1}{2}\tilde{G}_2(\tau)
\right)\tilde{y}(\tau).$$ If we iterate the previous formula we
obtain
$$\eta(\tau)\left(q \frac{d}{dq}\right)^2 y(\tau)=
\eta(\tau)\left(q \frac{d}{dq}\right) y'(\tau)$$
$$=\left(\left(q
\frac{d}{dq}\right)+\frac{1}{2}\tilde{G}_2(\tau) \right)
\eta(\tau)y'(\tau)=\left(\left(q
\frac{d}{dq}\right)+\frac{1}{2}\tilde{G}_2(\tau) \right)^2 y(\tau).$$
Now, by the induction, we get
$$\left(\left(q
\frac{d}{dq}\right)+\frac{1}{2}\tilde{G}_2(\tau) \right)^r
\tilde{y}(\tau)=\eta(\tau)y^{(r)}(\tau).$$ If we apply now the
Leibnitz rule we get \be \label{leibnitz}
\eta(\tau)y^{(r)}(\tau)=\left(q \frac{d}{dq} \right)^{r}
\tilde{y}(\tau)+\frac{r}{2} \tilde{G}_2(\tau) \left(q \frac{d}{dq}
\right)^{r-1} \tilde{y}(\tau)+\cdots, \ee where the dots denote
the terms with lower order derivatives of $\tilde{y}(\tau)$. The proof
now follows after we multiply the equation (\ref{diffm1}) by
$\eta(\tau)$ and apply (\ref{leibnitz}) for $r=1,...,k$. \epf
\begin{corollary}
We have  \be \label{wronskimod}
W(\tilde{y}_1(\tau),...,\tilde{y}_k(\tau))=\tilde{C}_k
\eta(\tau)^{(2k-1)k}, \ee where $\tilde{C}_k$ is some non-zero constant.
\end{corollary}
Now, let us work out the Wronskian on the left hand side of
equation (\ref{wronskimod}).
From (\ref{character}) it follows that
$$\tilde{y}_i(\tau)=q^{\left(h^{1,i}_{2,2k+1}-\frac{c_{2,2k+1}}{24}+\frac{1}{24}\right)}
\sum_{n_i \in \mathbb{Z}} (-1)^{n_i}
q^{\frac{(2(k-i)+1)n_i+(2k+1)n^2_i}{2}}.$$ For
every $n_i \in \mathbb{Z}$, $i=1,...,k$, let
$$a(n_i)=h^{1,i}_{2,2k+1}-\frac{c_{2,2k+1}}{24}+\frac{1}{24}
+\frac{(2(k-i)+1)n_i+(2k+1)n^2_i}{2},$$
$$\tilde{y}_i(\tau)=\sum_{n_i \in \mathbb{Z}} (-1)^{n_i} q^{a(n_i)}.$$
Now we prove Theorem \ref{theoremmain}. \\
\noindent {\em Proof of Theorem \ref{theoremmain}:} By
(\ref{wronskivander})
\bea &&
W(\tilde{y}_1(\tau),...,\tilde{y}_k(\tau))=
\left| \begin{array}{ccccc} \tilde{y}_1(\tau) & \tilde{y}_2(\tau) & . & . & \tilde{y}_k(\tau) \\
                              \tilde{y}'_1(\tau) & \tilde{y}'_2(\tau) & . & . & \tilde{y}'_k(\tau) \\
                              \tilde{y}''_1(\tau) & \tilde{y}''_2(\tau) & . & . &  \tilde{y}''_k(\tau) \\
                              . & . & . & . &  . \\
                              \tilde{y}^{(k-1)}_1(\tau) & \tilde{y}^{(k-1)}_2(\tau) & . & . & \tilde{y}^{(k-1)}_k(\tau)
\end{array} \right| \nn
&& =\sum_{n_1,...,n_k} \left| \begin{array}{ccccc}
(-1)^{n_1}q^{a(n_1)} & (-1)^{n_2}q^{a(n_2)} & . & .
&  (-1)^{n_k}q^{a(n_k)} \\
                              (-1)^{n_1}a(n_1)q^{a(n_1)} & (-1)^{n_2}a(n_2) q^{a(n_2)} & . & .
                              & (-1)^{n_k}a(n_k) q^{a(n_k)} \\
                              . & . & . & . &  . \\
                              . & . & . & . &  . \\
                              (-1)^{n_1}a(n_1)^{k-1} q^{a(n_1)} & (-1)^{n_2}a(n_2)^{k-1}
                              q^{a(n_2)} & . & . &
                              (-1)^{n_k}a(n_k)^{k-1}q^{a(n_k)}
\end{array} \right| \nn
&& =\sum_{n_1,...,n_k} (-1)^{\sum_{i=1}^k n_i} \left|
\begin{array}{ccccc} q^{a(n_1)} & q^{a(n_2)} &
. & .
&  q^{a(n_k)} \\
                              a(n_1)q^{a(n_1)} & a(n_2) q^{a(n_2)} & . & .
                              & a(n_k) q^{a(n_k)} \\
                              . & . & . & . &  . \\
                              . & . & . & . &  . \\
                              a(n_1)^{k-1} q^{a(n_1)} & a(n_2)^{k-1}
                              q^{a(n_2)} & . & . &
                              a(n_k)^{k-1}q^{a(n_k)}
\end{array} \right| \nn
&&=\sum_{n_1,...,n_k} (-1)^{\sum_{i=1}^k n_i} q^{\sum_{i=1}^k
a(n_i)} V(a(n_1),...,a(n_k)) \nn &&=(-1)^{\frac{k(k-1)}{2}}
\sum_{n_1,...,n_k} (-1)^{\sum_{i=1}^k n_i} q^{\sum_{i=1}^k a(n_i)}
\prod_{1 \leq i < j \leq k} (a(n_i)-a(n_j)). \eea Clearly,
$$\sum_{i=1}^k a(n_i)=\frac{2k^2-k}{24}+\sum_{i=1}^k \frac{(2(k-i)+1)n_i+(2k+1)n_i^2}{2}$$
and hence (cf. (\ref{bfn}))
$$L({\bf n})=\sum_{i=1}^k a(n_i).$$
Also
$$(-1)^{k(k-1)/2}\prod_{i < j}
(a(n_i)-a(n_j))=\prod_{i<j}(a(n_{k-i+1})-a(n_{k-j+1})).$$ Now the
formula
$$(a(n_{k-i+1})-a(n_{k-j+1}))$$
$$=\frac{1}{2(2k+1)}\left(i-j+ 2 k n_i-2k n_j +n_i-n_j \right)
\left(i+j-1+2n_i k+2 n_j k +n_i+n_j \right)$$ together with \bea
\label{chi} && \chi_D({\bf n})=\prod_{1 \leq i
<j \leq k}((2i-1+n_i(4k+2))^2-(2j-1+n_j(4k+2))^2) \\
&&= 2^{k(k-1)} \prod_{1 \leq i<j \leq k} (i-j+2n_i k -2 n_j
k+n_i-n_j)(1+j-1+2n_ik+2n_jk+n_i+n_j) \nonumber \eea
imply
$$\tilde{C}_k \eta(\tau)^{2k^2-k}=\frac{(4k+2)^{k(k-1)/2}}{2^{k(k-1)}}
q^{\frac{2k^2-k}{24}} \sum_{n_1,...,n_k} (-1)^{\sum_{i=1}^k n_i}
\chi_D({\bf n}) q^{L({\bf n})},$$
where the summation is over all $k$--tuples ${\bf n}=(n_1,...,n_k) \in \mathbb{Z}^k$. 
The constant $\tilde{C}_k$ is equal to
$$\frac{(4k+2)^{k(k-1)/2}}{\ds{\prod_{1 \leq i<j \leq k}(i-j)(i+j-1)}}$$
and this proves Theorem \ref{theoremmain}. \epf
\renewcommand{\theequation}{\thesection.\arabic{equation}}
\renewcommand{\thetheorem}{\thesection.\arabic{theorem}}
\setcounter{equation}{0} \setcounter{theorem}{0}

\section{General case}
In this part we extend specialized Macdonald's identities from the previous
section to all $c_{s,t}$--minimal models. 
(Un)fortunately, characters for minimal models in general
do not admit nice infinite product expansions. 
We shall rewrite first the character formula in Theorem \ref{struct} as a single sum (cf. \cite{HK}):
\be \label{genchar}
\bar{{\rm ch}}_{c_{s,t},h_{s,t}^{m,n}}(q)=\frac{\ds{\sum_{r \geq 0} \chi_{2st}^{h_{m,n}}(r) q^{\frac{r^2}{4 st}}} } {\eta(q)},
\ee
where
$$\chi_{2st}^{h_{m,n}}(r)=\left\{ \begin{array}{ccc} 1  \ {\rm for} \ 
r = \pm (ns - mt) \ {\rm mod} \ 2st \\
 - 1  \ {\rm for} \ r =  \pm (ns + mt) \ {\rm mod} \ 2st \\
0 , \ {\rm otherwise}. \end{array} \right. $$
Even though this formula is not very transparent it is the perhaps the only
closed expression that covers for all the minimal models.
Also, it seems that $\chi_{2st}^{h_{m,n}}$ does not have some obvious
arithmetic interpretation. 

Now, among all the pairs $(m,n)$ we have to identify those
pairs that give different $h_{s,t}^{m,n}$ values ( $\frac{(s-1)(t-1)}{2}$ in total).
It is not hard to see that the first
$$k=\frac{(s-1)(t-1)}{2}$$
values in the sequence 
$$h_{s,t}^{1,1}, h_{s,t}^{1,2},...,h_{s,t}^{s-1,t-2},h_{s,t}^{s-1,t-1},$$
starting with $h_{s,t}^{1,1}$, give the wanted values. For simplicity we will enumerate these 
(rational) numbers
by $$h_1,h_2,...,h_k.$$
We will need a stronger version of Theorem \ref{diffmain}.
\begin{theorem} \label{diffmainstrong}
Let $k$ be as above and for every $i=1,...,k$, let
$$\tilde{{y}}_i(\tau)=\eta(\tau) \bar{{\rm ch}}_{c_{s,t},h_i}(q),$$
Then
$$W(\tilde{y}_1(\tau),...,\tilde{y}_k(\tau))=C_{s,t} \eta(\tau)^{(2k-1)k},$$
where $C_{s,t}$ is some constant that depends only on $s$
and $t$.
\end{theorem}
{\em Proof:} In order to apply Theorem \ref{diffmain} for all
$c_{s,t}$--minimal models we need the following fact (observed by Feigin and Fuchs):
The vacuum module $V(c_{s,t},0)$ (see \cite{Wa}, \cite{M}), which carries a vertex
operator algebra structure, contains a singular vector of the
weight $(s-1)(t-1)$ of the form
$$(L^{\frac{(s-1)(t-1)}{2}}(-2)+ \cdots ){\bf 1} \in V(c_{s,t},0) $$
where the dots denote the lower order terms in the natural
filtration of $U({\rm Vir}_{\leq - 2})$ (for the notation see
\cite{M}, \cite{FFr}). Now, if apply the same procedure as in \cite{M}, we get 
a $k$--th order homogeneous linear differential equation 
satisfied by $\bar{{\rm ch}}_{c_{s,t},h_i}(q)$, $i=1,...,k$, which
is certainly a fundamental system of solutions (characaters for $c_{s,t}$--minimal models
are always linearly independent). This differential
equation is (again) of the form \be \label{gendiff}
\left(q\frac{d}{d \tau} \right)^{k}
y(\tau)+f_1(\tau)\left(q\frac{d}{d \tau} \right)^{k-1}y(\tau)+
\cdots+f_k(\tau) y(\tau)=0, \ee where
$$f_1(\tau)=k(k-1) \tilde{G}_2(\tau),$$
and $f_i(\tau)$, $i \geq 2$ are some polynomials in Eisenstein
series \cite{Zh}. Now we apply Lemma \ref{fixing} and the proof follows.
Let us notice here that existance of a singular vector in $V(c_{s,t},0)$ of the weight
$(s-1)(t-1)$ has been used in in \cite{DLM} (and assumed in \cite{Zh}) to show that
the vertex operator algebra $L(c_{s,t},0)$ satisfies the $C_2$-cofiniteness condition, which imples
regularity. \epf

\begin{theorem} ($c_{s,t}$--denominator formula) \label{theoremmain2}
For every $s,t, \in \mathbb{N}$, such that $2 \leq s < t$ and $(s,t)=1$, 
\be \label{ggg}
C_{s,t} \eta(\tau)^{2k^2-k}=\sum_{{\bf n} \in \mathbb{N}_0^k}
\left( \prod_{i=1}^k \chi_{2st}^{h_i} (n_i) \right)
V(n_1^2,...,n_k^2) q^{L({\bf n})},
\ee
where 
${\bf n}=(n_1,...,n_k) \in \mathbb{N}_0^k$,
 $$L({\bf
n})=\sum_{i=1}^k \frac{n_i^2}{4st}$$
and $C_{s,t}$ is some non-zero constant.
\end{theorem}
{\em Proof:} Theorem \ref{diffmainstrong} gives an
expression for $C_{s,t} \eta(\tau)^{2k(k-1)}$ as a Wronskian determinant. Therefore,
we only have to evaluate the Wronskian determinant and for this we
use (\ref{wronskivander}). Notice that we can factor the constant $\left(\frac{1}{4st} \right)^{k(k-1)/2}$
from the Vandermonde determinant $V(\frac{n_1^2}{4st},...,\frac{n_k^2}{4st})$.
The proof follows.
\epf

The constant $C_{s,t}$ can be computed explicitely by comparing the 
leading nonzero coefficients on both sides of (\ref{ggg}). 
In the special case $s=2$, $t=2k+1$, Theorem \ref{theoremmain2}
implies the formula in Theorem \ref{theoremmain}.
\begin{remark}
{\em Let $\mu(k)$ be the number of positive integer
solutions $s,t$ of the equation $$2k=(s-1)(t-1),$$
where $2 \leq s < t$, $s$ and $t$ are relatively prime. For instance $\mu(9)=3$, with the
solutions $s=2,t=19$, $s=3,t=10$ and $s=4,t=7$.
Our Theorem \ref{theoremmain2} implies that for every integer
$k \geq 2$ we have $\mu(k)$ expressions for
$$\eta(\tau)^{(2k-1)k}.$$ 
Is
there some reasonable explanation why the number $\mu(k)$ appears here?}
\end{remark}
\renewcommand{\theequation}{\thesection.\arabic{equation}}
\renewcommand{\thetheorem}{\thesection.\arabic{theorem}}
\setcounter{equation}{0} \setcounter{theorem}{0}

\section{Factorization of linear combinations of characters: $c_{3,4}$ example}
As we already noticed $c_{2,2k+1}$--minimal models are easier to
handle because the numerator formulas for characters are various
specialization of theta constants (so we can use the Jacobi
triple product identity \cite{A}). In general characters of minimal
models do not rise to ``nice'' 
infinite-product expressions (at least not in a
straightforward way).
In spite of this, a recent work of Bytsko and Fring (see \cite{BF} and
references therein) shows that certain {\em linear combinations} of
characters do have infinite product expansions (these factorization
properties are useful to prove certain Rogers-Ramanujan-type
formulas via dilogarithms \cite{BF}). Unfortunately, most of these
factorizations occur for some special modules among
$c_{s,t}$-minimal models so it is not clear what to expect in
general. 

From our point of view it is natural to consider linear combinations of characters 
because of the following elementary fact:
Let $T$ be an invertible matrix with entries being numbers, and
$$(h_1(\tau),...,h_k(\tau))=T (f_1(\tau),...,f_k(\tau)),$$
where $(\cdot,...,\cdot)$ denote a $k$-size column vector, then
\be \label{char}
W(h_1(\tau),...,h_k(\tau))=det(T) W (f_1(\tau),...,f_k(\tau)).
\ee
This simple observation and factorization properties of 
linear combinations can be used
to prove some nontrivial modular identities (see below).

Perhaps the most famous example of factorization of linear
combinations of characters occurs when the central charge is
$c_{3,4}=\frac{1}{2}$, i.e., the Ising model. This
model is unitary and well-understood. 
In order to connect these models with infinite products we 
recall the definition of (normalized) Weber's
functions \cite{We}:
$$\goth{f}(\tau)=q^{\frac{-1}{48}} \prod_{n \geq
0} (1 + q^{n + \frac{1}{2}}),$$
$$\goth{f}_1(\tau)=q^{\frac{-1}{48}} \prod_{n \geq
0} (1 - q^{n + \frac{1}{2}}),$$
$$\goth{f}_2(\tau)=q^{\frac{1}{24}} \prod_{n \geq 0}(1+q^{n+1}).$$
Our goal is to prove the following identity.
\begin{proposition} \label{weber}
We have
\be
256 \left| \begin{array}{ccc} \goth{f}(\tau) & \goth{f}_1(\tau) & \goth{f}_2(\tau) \\
\goth{f}'(\tau) & \goth{f}'_1(\tau) & \goth{f}'_2(\tau) \\
\goth{f}''(\tau) & \goth{f}''_1(\tau) & \goth{f}''_2(\tau)
\end{array} \right|=7 \eta(\tau)^{12}={7} \sqrt{\Delta}, \nonumber
\ee where $\Delta$ is the Ramanujan's discriminant function.
\end{proposition}
{\em Proof:} Our Theorem  \ref{diffmainstrong} (cf. formula
(\ref{gendiff})) implies that there is a $3$-rd order homogeneous
linear differential equation
$$\left(q \frac{d}{dq} \right)^3 y(q) + 6 \tilde{G}(q) \left(q \frac{d}{dq}
\right)^2 y(q) + f_2(q) \left(q \frac{d}{dq} \right) y(q)+f_3(q)y(q)=0, $$
with a fundamental system of solutions being
$$\bar{ch}_{1/2,0}(q), \ \bar{ch}_{1/2,1/2}(q),\ \ {\rm and} \
\ \bar{ch}_{1/2,1/16}(q).$$  On the other hand it is known (see
for instance \cite{KR}) that \bea &&
\bar{ch}_{c_{3,4},h_{3,4}^{1,2}}(q)=q^{\frac{1}{24}} \prod_{n \geq
0}(1+q^{n+1}), \nn && \bar{ch}_{c_{3,4},h_{3,4}^{1,1}}(q) \pm
\bar{ch}_{c_{3,4},h_{3,4}^{1,3}}(q)=q^{\frac{-1}{48}} \prod_{n
\geq 0} (1 \pm q^{n + \frac{1}{2}}). \nonumber \eea Because of (\ref{char}),
$$W(\goth{f}(\tau),\goth{f}_1(\tau),\goth{f}_2(\tau))=
C \
W(\bar{ch}_{1/2,0}(\tau),\bar{ch}_{1/2,1/2}(\tau),\bar{ch}_{1/2,1/16}(\tau)),$$
where $C$ is some non-zero constant. Now, Theorem \ref{diffmain}
implies that
$$W(\bar{ch}_{1/2,0}(\tau),\bar{ch}_{1/2,1/2}(\tau),\bar{ch}_{1/2,1/16}(\tau))$$
is a non-zero multiple of $\eta(\tau)^{12}=\sqrt{\Delta}$. By
comparing the leading coefficients in
$W(\goth{f}(\tau),\goth{f}_1(\tau),\goth{f}_2(\tau))$ and
$\eta(\tau)^{12}$ the proof follows. \epf
\begin{remark}
{\em We are confident that our Proposition \ref{weber} can
be proven by using classical Jacobi theta function techniques (e.g.
by using formulas for half-periods of Weierstrass series expressed
in terms of Jacobi theta constants and their derivatives).}
\end{remark}
\begin{remark}
{\em Our Proposition \ref{weber} implies that 
various factorization of linear combinations
(e.g., several formulas obtained in \cite{BF})
can be used to derive modular identities in terms of
Dedekind $\eta$-functions with different periods
(e.g. $\eta(m \tau )$, $\eta(\tau/n)$, $m,n \in \mathbb{N}$).}
\end{remark}
\renewcommand{\theequation}{\thesection.\arabic{equation}}
\renewcommand{\thetheorem}{\thesection.\arabic{theorem}}
\setcounter{equation}{0} \setcounter{theorem}{0}

\section{Conclusion and future work}
\begin{itemize}
\item[(i)] We obtained a two-parametric generalization of the Euler's identity (\ref{euler}).
In the case of $c_{2,2k+1}$--minimal models we gave a new proof of 
a series of specialized Macdonald's identities for the affine root system 
of type $A_{2k}^{(2)}$, $k \geq 2$, involving powers of the Dedekind
$\eta$--function (cf. Theorem \ref{theoremmain}). We stress that  
for different values of $s$ and $t$ we got
different looking identities. It is an open question to relate
identities in Theorem \ref{theoremmain2} to other specialized Macdonald's identities.
\item[(ii)] 
(Work in progress) We are modifying some of techniques from this
paper and \cite{M1}, \cite{M} for $N=1$ and $N=2$ superconformal minimal 
models in connection with a work of Kac and Wakimoto \cite{KW2} on affine Lie superalgebras, 
and Milne's work on sums of squares \cite{Mi}.
\end{itemize}
\small{

}
\vskip 5mm \noindent {\em \small \sc Department of Mathematical
Sciences, Rensselaer Polytechnic Institute, Troy, NY 12180 } \\
\noindent {\em E-mail address}: milasa@rpi.edu

\end{document}